\documentclass[12pt]{article}

\usepackage[english]{babel}
\usepackage{graphicx}
\usepackage{mathtools}
\usepackage[toc,page]{appendix}
\usepackage{amsmath,amsfonts,latexsym}

\usepackage{pdfpages}
\usepackage{amssymb}
\usepackage{graphics}
\usepackage{geometry}
\usepackage{lipsum}     
\usepackage[toc,page]{appendix}
\usepackage{fancyhdr}
\usepackage{setspace}

\newcommand{\nc}{\newcommand}
\nc{\tr}{{\vartriangle}} \nc{\vth}{{\vartheta}}
\nc{\bt}{{\beta}}
\nc{\dl}{{\delta}} \nc{\Dl}{{\tr}}
\nc{\p}{{\psi}}
\nc{\gm}{{\gamma}} \nc{\Gm}{{\Gamma}} \nc{\sg}{{\sigma}}
\nc{\ve}{{\varepsilon}} 
\nc{\ch}{{\cal H}} 
\nc{\cf}{{\cal F}}
\nc{\cp}{{\cal P}}
\nc{\td}{\tilde}
\nc{\ck}{\cal K}

\newtheorem{lemma}{Lemma}[section]

\newtheorem{corollary}[lemma]{Corollary}
\newtheorem{proposition}[lemma]{Proposition}

\newtheorem{example}[lemma]{Example}

\newtheorem{remark}[lemma]{Remark}
\numberwithin{equation}{section}

\begin{document}

\title
{Rate of Convergence of Truncated Stochastic Approximation Procedures with Moving Bounds}
   \author{Teo Sharia and Lei Zhong}
\date{}
\maketitle
\begin{center}
{\it \footnotesize
Department of Mathematics, Royal Holloway,  University of London\\
Egham, Surrey TW20 0EX \\ e-mail: t.sharia@rhul.ac.uk }
\end{center}


\begin{abstract}
The paper is concerned  with  stochastic approximation procedures 
having  three main characteristics:  truncations  with random moving bounds, 
 a matrix valued random step-size sequence,  and  a dynamically changing random regression function.
 We study convergence and rate of convergence. 
 Main results are supplemented with corollaries to establish various sets of sufficient conditions, with the  main emphases  on  the  parametric statistical estimation.  The  theory  is illustrated by examples and special cases.
 \end{abstract}


\begin{center}
\end{center}

\section{Introduction}\label{Intro}

This paper is a continuation of Sharia (2014) where a large class of truncated Stochastic approximation (SA) procedures with moving random bounds was proposed. Although the proposed class of procedures can be  applied to  a wider range of problems, our main motivation comes from applications to  parametric statistical estimation theory. To make this paper self contained, we introduce the main ideas below (a full list of references as well as some comparisons can be found in Sharia (2014)).

The main idea can be easily explained in the case of the  classical problem of finding a unique zero, say $z^0$,  of  a real valued function $R(z):  \mathbb{R}  \to \mathbb{R}$ when only noisy measurements  of $R$ are available. To estimate $z^0$, consider a sequence defined recursively   as                         
\begin{equation}\label{S}
Z_t=
Z_{t-1}+\gm_t \left[R(Z_{t-1})+\ve_t\right],   \qquad t=1,2,\dots
\end{equation}
where $\{\ve_t\}$ is a sequence of zero-mean random variables and  $\{\gamma_t\}$ is a deterministic sequence of positive numbers.
This is the classical Robbins-Monro SA procedure (see  Robbins and Monro (1951)\nocite{RM}), which under certain conditions  converges to the root $z^0$ of the equation $R(z)=0$.  (Comprehensive surveys of the SA technique can be found in Benveniste et al. (1990)\nocite{Ben1990}, Borkar (2008)\nocite{Bor},
Kushner and Yin (2003)\nocite{KushYin}, Lai (2003)\nocite{Lai}, and Kushner (2010)\nocite{Kush1}.)

Statistical parameter estimation is one of the most important applications of the above procedure. Indeed, suppose that $X_1, \dots, X_t$    are    i.i.d. random variables and $f(x, \theta)$ is the  common probability density  function (w.r.t. some $\sigma$-finite measure), where $\theta \in \mathbb{R}^m $ is an unknown parameter. Consider a recursive estimation procedure for $\theta$ defined by
\begin{equation} \label{reciid}
\hat\theta_t=\hat\theta_{t-1}+
\frac  {1}{t} i(\hat\theta_{t-1})^{-1}~\frac{{f'}^T(X_t, \hat\theta_{t-1})}{f(X_t, \hat\theta_{t-1})},
~~~~~~~~~t\ge 1,
\end{equation}
 where $\hat\theta_0\in {\mathbb{R}}^m$ is some starting value and $i(\theta)$ is the one-step Fisher information matrix ($f'$ is the row-vector of partial derivatives of $f$ w.r.t. the components of $\theta$). This estimator was introduced in Sakrison (1965\nocite{Sakr}) and studied by a number of authors (see e.g, Polyak and Tsypkin (1980)\nocite{Pol}, Campbell (1982)\nocite{Cam},  Ljung and Soderstrom (1987)\nocite{LaS},  Lazrieve and Toronjadze (1987)\nocite {lazTor},   Englund et al (1989)\nocite{Eng},  Lazrieve et al (1997, 2008)\nocite{Laz2}\nocite{Laz}, 
 Sharia (1997--2010)).\nocite{Shar0}\nocite{Shar}\nocite{Shar1} \nocite{Shar2}\nocite{Shar3}\nocite{Shar5} In particular, it has been shown
that under certain conditions, the recursive estimator $\hat\theta_t$ is asymptotically equivalent to the maximum likelihood estimator, i.e., it is  consistent and asymptotically efficient.  
One can analyse \eqref{reciid}
by rewriting it in the form of stochastic approximation with $\gamma_t=1/t$,
$$
R(z)=i(z)^{-1} E_{\theta}\left\{ \frac{{f'}^T(X_t, z)}{f(X_t, z)}\right\}
~~~ \mbox{and} ~~~ 
\ve_t=i(\hat\theta_{t-1})^{-1} \left(\frac{{f'}^T(X_t, \hat\theta_{t-1})}{f(X_t, \hat\theta_{t-1})} - R(\hat\theta_{t-1})\right),
$$ 
where $\theta$ is an arbitrary but fixed value of the unknown parameter. 
Indeed, under certain standard assumptions, $R(\theta)=0$ and $\{\ve_t\}$ is a martingale difference w.r.t. the filtration $\{{\cf}_{t}\}$ generated by  $\{X_t\}$. So, \eqref{reciid} is a standard SA of type \eqref{S}.

Suppose now that we have a stochastic process $X_1,  X_2, \dots $ and let 
$f_t(x,\theta)=f_t(x,\theta|X_1, \dots,X_{t-1})$ be  the conditional probability density function 
of the observation $X_t$ given $X_1, \dots,X_{t-1}$,  where $\theta \in \mathbb{R}^m $ is an unknown parameter.
 Then one can define a recursive estimator of   $\theta$ by
\begin{equation} \label{rec11}
\hat\theta_t=\hat\theta_{t-1}+
{\gamma_t(\hat\theta_{t-1})}\psi_t(\hat\theta_{t-1}),
~~~~~~~~~t\ge 1,
\end{equation}
where   $\psi_t(\theta)=\psi_t(X_1, \dots, X_t; \theta),$  $t=1,2,\dots, $
are suitably chosen functions which may, in general, depend on the  vector of all past and present observations $X_1,..., X_t$, 
and have the  property that   the process $\psi_t(\theta)$ is  $P^\theta$- martingale difference, i.e.,
$E_\theta\left\{ \psi_t(\theta)\mid{\cal{F}}_{t-1}\right\}=0$ for each $t$. For example,  a choice
$$\psi_t(\theta)=l_t(\theta)\equiv \frac{[f'_t(X_t,\theta)]^T}{f_t(X_t,\theta)} $$
  yields a   likelihood type estimation procedure. 
In general,   to obtain   an  estimator with asymptotically optimal properties,  a  state-dependent matrix-valued random step-size sequences are needed (see Sharia (2010\nocite{Shar3})). 
For  the above procedure,  a step-size sequence $\gamma_t(\theta)$
with the property
$$
\gamma_t^{-1}(\theta)-\gamma_{t-1}^{-1}(\theta)=E_{\theta}\{ \p_t(\theta)l^{T}_t(\theta)\mid
{\cf}_{t-1}\}
$$  
is an optimal choice. For example, to derive a recursive procedure which is asymptotically equivalent to the maximum likelihood  estimator,  we need to take  
$$\psi_t(\theta)=l_t(\theta)
 \;\;\;\mbox{ and }\;\;\;
\gamma_t(\theta)=I_t^{-1}(\theta),
$$
where   
\begin{equation} \label{fish}
I_t(\theta)=\sum_{s=1}^t E\{l_s(\theta)l_s^T(\theta)|{\cal F}_{s-1}\}
\end{equation}
is the conditional Fisher information  matrix. To rewrite \eqref{rec11} in the SA form, let us assume that $\theta$ is an arbitrary but fixed value of the parameter and define
 $$
R_t(z)=E_{\theta}\left\{ \psi_t(X_t, z)\mid{\cf}_{t-1}\right\}
~~~ \mbox{and} ~~~ 
\ve_t(z)=\left(\psi_t(X_t, z) - R_t(z)\right).
$$  
Then, since $\psi_t(\theta)$ is $P^\theta$-martingale difference, it follows that $R_t(\theta)=0$  for each $t$. So, the objective now is to find a common root $\theta$ of a dynamically changing sequence of functions $R_t$.

Before introducing the general SA process, let us consider one simple modification of the classical SA procedure. Suppose that we have additional information about the root $z^0$ of the equation $R(z)=0$. Let us, e.g., assume that  $z^0\in[\alpha_t,\beta_t]$ at each step $t$, where $\alpha_t$ and $\beta_t$ are random variables such that $-\infty<\alpha_t\leq\beta_t<\infty$. Then one can consider a procedure, which at each step $t$ produces points from the interval $[\alpha_t,\beta_t]$. For example, a truncated classical SA procedure in this case can be derived using the following recursion
\begin{equation}\label{STr}
Z_t=
\Phi_{[\alpha_t, \beta_t]}\big(~
Z_{t-1}+\gm_t \left[R(Z_{t-1})+\ve_t\right]\big),   \qquad t=1,2,\dots \nonumber
\end{equation}
where $\Phi$ is the truncation
operator, that is, for any $ -\infty<a \leq b<\infty ,$  ~
$$
\Phi_{[a,b]}(z) =
\begin{cases}
   a & \text{if} \;\; z<a,  \\
    z& \text{if} \;\;  a\le z\le  b,\\
    b & \text{if} \;\;z>b .
\end{cases}
$$

Truncated procedures  may be useful  in a number of circumstances. For example,
 if the functions in the recursive equation are defined only for  certain values of the parameter, then  
the procedure should  produce points only from this set. Truncations  may  also be useful when certain  standard assumptions, e.g., 
  conditions on  the  growth rate of the relevant functions are not satisfied.   
Truncations may  also help to make an efficient use of   auxiliary information  concerning the value of the unknown parameter.
For example, we might have auxiliary information about the parameters, e.g. a set, possibly time dependent, that  contains the value of the unknown parameter. Also, sometimes 
 a consistent but not necessarily efficient auxiliary estimator $\tilde\theta_t$  is available having a  rate $d_t$. Then
  to obtain asymptotically efficient estimator, one can construct a  procedure with  shrinking  bounds  by 
truncating  the recursive procedure in a neighbourhood of $\theta$
 with $ [\alpha_t, \beta_t]=[\tilde\theta_t-\delta_t, \tilde\theta_t+\delta_t],$ where $\delta_t \to 0$.  

Note that the idea of truncations is not new and goes back to Khasʹminskii and Nevelson (1972)\nocite{Khas} and  Fabian (1978)\nocite{Fab}  (see also 
 Chen and Zhu (1986)\nocite{Chen2}, Chen et al.(1987)\nocite{Chen1},  Andrad{\'o}ttir (1995)\nocite{Andr}, Sharia (1997)\nocite{Shar0},  Tadic (1997,1998)\nocite{Tadic1}\nocite{Tadic2}, Lelong (2008)\nocite{Lel}.
A comprehensive bibliography and some comparisons can be found in Sharia (2014)\nocite{Shar4}).

In order to study these procedures in an  unified  manner,  Sharia (2014)\nocite{Shar4} introduced  a SA of the following form 
$$                         
Z_t=\Phi_{U_t} \Big(~
Z_{t-1}+\gm_t(Z_{t-1}) \big[ R_t(Z_{t-1})+\ve_t(Z_{t-1}) \big]\Big),   \quad t=1,2,\dots
$$
where  $Z_0 \in \mathbb{R}^m$  is some starting value,  $R_t(z)$ is  a predictable process with the property that $R_t(z^0)=0$ for all $t$'s, $\gm_t(z)$ is a matrix-valued predictable step-size sequence, and $U_t   \subset \mathbb{R}^m$  is  a random sequence of truncation sets (see Section \ref{MON} for details).
These SA procedures have the  following main characteristics: (1) inhomogeneous  random functions $R_t$; (2) state dependent matrix valued random step-sizes;  (3)  truncations with random and  moving 
(shrinking or expanding) bounds.
The main motivation for these comes from parametric statistical 
applications: (1) is needed for recursive parameter estimation procedures for  non i.i.d. models;
(2) is required to guarantee asymptotic optimality and efficiency of  statistical estimation; (3) is needed for
various different adaptive truncations, in particular, for the ones arising by auxiliary estimators. 

 Convergence of the above class of procedures is studied in Sharia (2014). In this paper we present  new results on  rate of convergence.  Furthermore, we  present a convergence result which generalises the corresponding result in Sharia (2014)\nocite{Shar4} by considering time dependent random Lyapunov type functions (see Lemma \ref{RC}). This generalisation turns out to be quite useful as it can be used to derive convergence results of the recursive parameter estimators in  time series models.   Some of the conditions in the main statements are difficult to interpret. Therefore, we discuss these conditions in  explanatory remarks and corollaries. The corollaries are presented in such a way that each subsequent statement imposes conditions that are more restrictive than the previous one.
We discuss  the case of the classical SA and  demonstrate that conditions introduced in this paper are minimal in the sense that  they do not impose  any additional restrictions when applied to  the classical case. 
We also compare our set of conditions to that of Kushner-Clark's setting (see Remark \ref{K-C}).
Furthermore, the paper contains new results even for the classical SA. In particular, truncations with moving bounds give a possibility to use SA in the cases when the standard conditions on the function $R$ do not hold.
Also, an interesting link between the rate of the step-size sequence and the rate of convergence of the SA process is given  in the classical case (see corollary \ref{Po} and Remark \ref{Po1}). This observation might not surprise experts working in this field, but we failed to find it in a written form in the existing literature.

\section{Main objects and notation}\label{MON}

Let $(\Omega, ~ \cf,F=(\cf_t)_{t\geq 0}, ~P)$ be a stochastic basis satisfying the usual conditions. Suppose that 
for each $t=1,2, \dots$, we have 
$( {\cal{B}}  ( \mathbb{R}^m)  \times \cf )$-measurable functions
$$
\begin{array}{cl}
 R_t(z)= R_t(z,\omega) &:\mathbb{R}^m \times  \Omega   \to    \mathbb{R}^m  \\
 \ve_t(z)=\ve_t(z,\omega) &:\mathbb{R}^m \times  \Omega   \to    \mathbb{R}^m      \\
  \gamma_t(z)=\gamma_t(z,\omega)&:\mathbb{R}^m \times  \Omega   \to    \mathbb{R}^{m\times m}    
\end{array}
$$
such that 
 for each $z\in  \mathbb{R}^m$,  the   processes    $R_t(z) $  and  $\gamma_t(z)$ are  predictable, i.e.,
 $R_t(z) $  and  $\gamma_t(z)$ are $\cf_{t-1}$ measurable for each $t$.
 Suppose also that  
 for each $z\in  \mathbb{R}^m$,  the process $\ve_t(z) $   is a martingale difference, i.e., $\ve_t(z) $ 
 is  $\cf_{t}$ measurable and  $E\left\{\ve_t(z)\mid{\cal{F}}_{t-1}\right\}=0$.
 We also assume that 
$$
R_t(z^0)=0
$$  
for each $ t=1, 2, \dots $, where  $z^0 \in   \mathbb{R}^m$   is  a non-random vector. 

Suppose that $h=h(z)$ is a real valued function  of
 $z \in {{\mathbb{R}}}^m$.  Denote by $ h'(z)$  the row-vector
 of partial derivatives
of $h$ with respect to the components of $z$, that
is,
 $
  h'(z)=\left(\frac{{\partial}}{{\partial} z_1} h(z), \dots,
 \frac{{\partial}}{{\partial} z_m} h(z)\right).
 $
Also, we denote by  $h''(z)$ the  matrix of second partial derivatives.
 The $m\times m$ identity matrix is denoted by ${{\bf I}}$.
Denote by $[a]^+$ and $[a]^-$ the positive and negative parts of $a\in \mathbb R$, i.e. $[a]^+=\max(a,0)$ and $[a]^-=\min(a,0)$.

Let    $U \subset    \mathbb{R}^m$ is a closed convex set  and define a truncation  operator    as a function
$\Phi_U(z) : \mathbb{R}^m  \longrightarrow \mathbb{R}^m$, such that
$$
\Phi_U(z) =\begin{cases}
    z & \text{if} \;\; z\in U \\
     z^* & \text{if} \;\; z\notin U,
\end{cases}
$$
where $z^*$ is a point in 
$U$,  that minimizes   the distance  to $z$.

Suppose that   $z^0 \in   \mathbb{R}^m$.  We say that a random sequence   of sets $U_t =U_t(\omega)$ 
($ t=1,2, \dots $)  from $\mathbb{R}^m $  is {\underline{\bf admissible}}  for  $z^0$  if
\medskip

\noindent
$\bullet$ 
for each $t$ and  $\omega,$    $U_t(\omega)$   is a   closed convex  subset  of $ \mathbb{R}^m$;
 \\
$\bullet$ 
 for each   $t$ and $z \in   \mathbb{R}^m$, the truncation $\Phi_{U_t}(z) $  is  $ {\cal{F}}_{t}$ measurable; 
 \\
$\bullet$   $z^0\in U_t$ eventually, i.e.,  for almost all  $\omega$ there exist  $t_0(\omega)<\infty$
 such that $z^0\in U_t(\omega)$ whenever $t >t_0(\omega)$.

\medskip
Assume that $Z_0 \in \mathbb{R}^m$  is some starting value and  consider the procedure 
\begin{equation}\label{TSA}
Z_t=
\Phi_{U_t}
\Big(
Z_{t-1}+\gm_t(Z_{t-1}) \Psi_t(Z_{t-1})\Big),   \quad t=1,2,\dots
\end{equation}
where  $U_t $   is {admissible}   for  $z^0$,
$$
\Psi_t(z)=R_t(z)+\ve_t(z),
$$
and   $R_t(z) $, $\ve_t(z)$, $\gm_t(z)$ are random fields  defined above. Everywhere in this work, we assume that
\begin{equation}\label{GTSA1}
 E\left\{\Psi_t(Z_{t-1})\mid{\cal{F}}_{t-1}\right\}=R_t(Z_{t-1})
  \end{equation}
                     and
\begin{equation}\label{GTSA2}
 E\left\{\ve_t^T(Z_{t-1})\ve_t(Z_{t-1})\mid{\cal{F}}_{t-1}\right\}= 
 \left[E\left\{\ve_t^T(z)\ve_t(z)\mid{\cal{F}}_{t-1}\right\} \right] _{z=Z_{t-1}},
 \end{equation}
 and the conditional expectations   \eqref{GTSA1} and \eqref{GTSA2}  are assumed to be finite.
 \medskip

%
\begin{remark}\label{disint}  {\rm Condition \eqref{GTSA1} ensures that $\ve_t(Z_{t-1})$ is a martingale difference.
Conditions  \eqref{GTSA1} and \eqref{GTSA2} obviously hold if, e.g.,  the  measurement errors $\ve_t(u)$ are independent 
random variables, or if they  are state independent. In general,  
since we assume that all  conditional  expectations are calculated  as integrals w.r.t. corresponding regular conditional probability measures (see the convention below), these conditions can be checked using disintegration formula (see, e.g.,  Theorem 5.4 in Kallenberg (2002)\nocite{Kall}).}
\end{remark}

We say that a random field
$$
V_t(z)=V_t(z,\omega): {\mathbb R}^m\times \Omega\longrightarrow {\mathbb R}\;\;\;\;\;\;(t=1,2,...)
$$
is a {\underline{\bf Lyapunov random field}} if
\\
\\
\noindent$\bullet$ $V_t(z)$ is a predictable process for each $z\in{\mathbb R}^m$;

\noindent$\bullet$ for each $t$ and almost all $\omega$, $V_t(z)$ is a non-negative function with continuous and bounded partial second derivatives.
\\

\noindent
{\bf \em Convention.}

\noindent
$\bullet$
 {\em Everywhere in the present work
  convergence and all relations between random
variables are meant with probability one w.r.t. the measure
$P$ unless specified otherwise. \\
{$\bullet$} A sequence of random
variables $(\zeta_t)_{t\ge1}$ has a property  {\underline {{\bf \em eventually}}} if for
every $\omega$ in a set $\Omega_0$ of $P$ probability 1, the realisation 
 $\zeta_t(\omega)$  has this property for all $t$ greater than some
$t_0(\omega)<\infty$.}\\
{$\bullet$} {\em All conditional expectations are calculated  as integrals w.r.t. corresponding regular conditional probability measures.}\\
{$\bullet$}
{\em The $\inf_{z\in U} h(z)$ of a real valued function $h(z)$ is  $1$
whenever $U=\emptyset$.}%


 \section{Convergence and rate of convergence}\label{MR}

We start this section with a convergence lemma, which uses a concept of a Lyapunov random field (see Section \ref{MON}). The proof of this lemma is very similar to that of presented in Sharia (2014). However, the dynamically changing Lyapunov functions make it possible to apply this result to derive the rate of convergence of the SA procedures. Also, this result turns out to be very useful to derive convergence of the recursive parameter estimations in time series models.

 \begin{lemma}\label{RC} 
 Suppose that  $Z_t$  is a  process  defined by \eqref{TSA}. Let $V_t(u)$ be a Lyapunov random field. Denote $\Delta_t=Z_t-z^0$, $\Delta V_t(u)=V_t(u)-V_{t-1}(u)$, and assume that
\begin{description}
\item[(V1)]\label{V1}
$$
V_t(\Delta_t)\leq V_t\Big( \Delta_{t-1}+\gamma_t(Z_{t-1})[R_t(Z_{t-1})+\ve_t (Z_{t-1}) ]\Big)
$$
eventually;
\item[(V2)]\label{V2}
\begin{equation*}\label{CO}
\sum_{t=1}^\infty [1+V_{t-1}(\Delta_{t-1})]^{-1} [{\ck}_t(\Delta_{t-1})]^+<\infty,\qquad
 P\mbox{-a.s.},
\end{equation*}
where
$$
{\ck}_t(u)= \Delta V_t(u)+V_t'(u) \gm_t(z^0+u)R_t(z^0+u)+\eta_t(z^0+u)
$$
and
$$
\eta_t(v)=\frac 1 2 \sup_z E
\left\{\Big[R_t(v)+{\ve}_t(v)\Big]^T\gm_t^T(v) V_t''(z)\gm_t(v) \Big[R_t(v)+{\ve}_t(v)\Big] \Big|
{\cf}_{t-1}\right\}. \nonumber
$$
\end{description}
Then $V_t(\Delta_t)$ converges  ($P$-a.s.) to a finite limit for any initial value $Z_0$.
\\
\\
Furthermore, if there exists a set A $\in {\cf}$ with $P(A)>0$ such that for each $\epsilon\in (0,1)$
\begin{description}
\item[(V3)]
\begin{equation}\label{CO2}
\sum_{t=1}^\infty \inf_{\stackrel{\epsilon\leq V_t(u)\leq1/\epsilon}{z^0+u\in U_{t-1}}}[{\ck}_t(u)]^-=\infty \quad \mbox{on A},
\end{equation}
\end{description}
then $V_t(\Delta_t)\longrightarrow 0$ ($P$-a.s.) for any initial value $Z_0$.
\end{lemma}
{\bf Proof.} The proof is similar to that of Theorem 2.2 and 2.4 in Sharia (2014)\nocite{Shar4}. Rewrite \eqref{TSA} in the form
$$
\Delta_t=\Delta_{t-1}+{\gm}_t(Z_{t-1})[R_t(Z_{t-1})+{\ve}_t(Z_{t-1})].
$$
By (V1), using the Taylor expansion, we have
\begin{eqnarray}
V_t(\Delta_t)&\leq&V_t(\Delta_{t-1})+V'_t(\Delta_{t-1}){\gm}_t(Z_{t-1})[R_t(Z_{t-1})+{\ve}_t(Z_{t-1})]\nonumber\\
&&
+\frac 1 2 [R_t(Z_{t-1})+{\ve}_t(Z_{t-1})]^T {\gm}_t^T(Z_{t-1}) V''_t(\tilde{ \Delta}_{t-1}){\gm}_t(Z_{t-1})[R_t(Z_{t-1})+{\ve}_t(Z_{t-1})],\nonumber
\end{eqnarray}
where $\tilde{\Delta}_{t-1}\in \mathbb R^m $ is ${\cal F}_{t-1}$-measurable
Since
$$
V_{t}(\Delta_{t-1})=V_{t-1}(\Delta_{t-1})+\Delta V_t(\Delta_{t-1}),
$$
using \eqref{GTSA1} and \eqref{GTSA2}, we obtain
$$
E\{V_t(\Delta_t)|{\cf}_{t-1}\}\leq V_{t-1}(\Delta_{t-1})+{\ck}_t(\Delta_{t-1}).
$$
Then, using the decomposition ${\ck}_t=[{\ck}_t]^+-[{\ck}_t]^-$, the above can be rewritten as
$$
E\{V_t(\Delta_t)|{\cf}_{t-1}\}\leq V_{t-1}(\Delta_{t-1})(1+B_t)+B_t-[{\ck}_t(\Delta_{t-1})]^-,
$$
where $B_t=(1+V_{t-1}(\Delta_{t-1}))^{-1}[{\ck}_t(\Delta_{t-1})]^+$.

By $(V2)$, we have that $\sum_{t=1}^\infty B_t <\infty$. Now we can use Lemma \ref{RMTh} in Appendix  (with $X_t=V_t(\Delta_t), \beta_{t-1}=\xi_{t-1}=B_t$ and $\zeta_t=[{\ck}_t(\Delta_{t-1})]^-$) to deduce that the processes $V_t(\Delta_t)$ and
$$
Y_t=\sum_{s=1}^t [{\ck}_s(\Delta_{s-1})]^-
$$
converge to some finite limits.
Therefore, it follows that $V_t(\Delta_t) \rightarrow r \geq 0$. 

To prove the second assertion, suppose that $r>0$. Then there exist $\epsilon>0$ such that $\epsilon\leq V_t(\Delta_t)\leq 1/\epsilon$ eventually. By \eqref{CO2}, this would imply that for some $t_0$,
$$
\sum_{s=t_0}^\infty [{\ck}_s(\Delta_{s-1})]^-\geq \sum_{s=t_0}^\infty \inf_{\stackrel{\epsilon\leq V_s(u)\leq1/\epsilon}{z^0+u\in U_{s-1}}}[{\ck}_s(u)]^-=\infty
$$
on the set A, which contradicts the existence of a finite limit of $Y_t$. Hence, $r=0$ and $V_t(\Delta_t) \longrightarrow 0$. \hfill $\blacksquare$
\begin{remark}{\rm
The conditions of the above Lemma are difficult to interpret. Therefore, the rest of the section is devoted to formulate lemmas and corollaries
(Lemmas \ref{CL} and \ref{EGV}, Corollaries \ref{SC}, \ref{EGV1D} and \ref{Dev})  containing sufficient conditions for the convergence and the rate of convergence, 
and remarks  (Remarks \ref{FTC0}, \ref{tsphere}, \ref{FTC}, \ref{FTC2}, \ref{FTC3} and \ref{FTC4}) explaining  some of the assumptions.  These results are presented in such a way, that each subsequent statement imposes conditions that are more restrictive than the previous one. For example, Corollary \ref{Dev} and Remark \ref{FTC4} contain conditions which are most restrictive than all the previous ones, but are written in the simplest possible terms.
}\end{remark}
\begin{remark}\label{FTC0}{\rm
A typical choice of $V_t(u)$ is $V_t(u)=u^T C_t u$, where $\{C_t\}$ is a  predictable positive semi-definite matrix process. If $C_t/a_t$ goes to a finite matrix with $a_t\longrightarrow \infty$, then subject to the conditions of Lemma \ref{RC}, $a_t\|Z_t-z^0\|^2$ will tend to a finite limit implying that $Z_t\longrightarrow z^0$. This approach is adopted in Example \ref{RLS} to derive convergence of the on-line Least Square estimator.
}\end{remark}
\begin{remark}\label{tsphere}{\rm
Consider truncation sets $U_t=S(\alpha_t, r_t)$, where $S$ denotes a closed sphere in $\mathbb R^m$ with the center at $\alpha_t\in\mathbb R^m$ and the radius $r_t$. Let $z_t'=\Phi_{U_t}(z_t)$ and suppose that $z^0\in U_t$. Let $V_t(u)=u^TC_tu$ where $C_t$ is a positive definite matrix and denote by $\lambda_t^{max}$ and $\lambda_t^{min}$ the largest and smallest eigenvalues of $C_t$ respectively. Then $(z_t'-z^0)^TC_t(z_t'-z^0)\leq(z_t-z^0)^TC_t(z_t-z^0)$  \Big(i.e., (V1) holds with $V_t(u)=u^TC_tu$\Big), if $\lambda_t^{max}v_t^2\leq\lambda_t^{min} r_t^2$, where $v_t=\|\alpha_t-z^0\|$. (See Proposition \ref{TsphPr} in Appendix  for details.) In particular, if $C_t$ is a scalar matrix, condition (V1) automatically holds.
}\end{remark}

\begin{lemma} \label{CL} 
Suppose that all the conditions of Lemma \ref{RC} hold and
\begin{description}
\item[(L)]for any $M>0$, there exist some $\delta=\delta(\omega)>0$ such that 
\begin{equation*}\label{VV}
\inf_{ \|u\| \geq M} V_t(u)> \delta\;\;\;\;\;\; \mbox{ eventually.}
\end{equation*}
\end{description}
Then   $Z_t \longrightarrow z^0 \;\;
(P$-a.s.)  for any  initial value $Z_0$.
\end{lemma}
{\bf Proof.} From Lemma \ref{RC}, we have $V_t(\Delta_t)\longrightarrow 0$ (a.s.). Now, $\Delta_{t} \longrightarrow  0$ follows from (L) by contradiction.  Indeed, suppose that  $\Delta_{t} \not\longrightarrow  0$ on a set, say $B$ of positive probability. Then, for any fixed $\omega$ from this set, there would exist a sequence  $t_k \longrightarrow \infty$  such that $\|\Delta_{t_k}\| \ge \epsilon $ for some $\epsilon>0,$ and \eqref{VV} would imply that $ V_{t_k}(\Delta_{t_k}) >  \delta>0$ for large $k$-s, which contradicts the  $P$-a.s.  convergence   $ V_t(\Delta_{t}) \longrightarrow  0$.\hfill $ \blacksquare $


  \begin{remark}\label{PoD}{\rm
  The following corollary contains  simple sufficient conditions for convergence. The poof of this  corollary does not require dynamically changing Lyapunov functions and can be obtained from a less general version of Lemma \ref{RC} presented in Sharia (2014). We decided to present this corollary for the sake of completeness, noting that the proof, as well as 
 a number of different sets of sufficient conditions, can be found   in Sharia (2014\nocite{Shar4}).  
 } \end{remark}

\begin{corollary}\label{SC}   Suppose that  $Z_t$  is a  process  defined by 
\eqref{TSA}, $U_t$ are admissible truncations for $z^0$ and 

\begin{description} 
 \item[(D1)] for large $t$'s 
 $$
 (z-z^0)^T R_t(z) \le 0 \;\;\;\mbox{if}\;\;\;z \in U_{t-1};
 $$
\item[(D2)]there exists a predictable process $r_t>0$ such that
$$
  \sup_{z \in U_{t-1}}\frac {E \left\{ \|R_t(z)+\ve_t(z)\|^2 \mid{\cf}_{t-1}\right\}}
  {1+\| z-z^o\|^2}\le r_t $$    
 eventually, and
 $$
  \sum_{t=1}^{\infty}   {r_{t}}{a_t^{-2}} <\infty, \qquad P\mbox{-a.s.}
$$
 \end{description}
Then $\|Z_t-z^0\|$ converges  ($P$-a.s.) to a finite limit.

\bigskip
\noindent Furthermore, if 
\begin{description}
 \item[(D3)]   for each $\epsilon\in (0, 1),$ there exists a predictable process $\nu_t>0$ such that
 \begin{equation} \label{INF}
 \inf_{\stackrel{ \epsilon \le \|z-z^o\| \le 1/\epsilon}{z\in U_{t-1}}}  -(z-z^0)^T R_t(z)> \nu_t \nonumber
\end{equation}
eventually, where
$$
  \sum_{t=1}^{\infty}  {\nu_{t}}{a_t^{-1}} =\infty, \qquad P\mbox{-a.s.}
$$
\end{description}
Then $Z_t$ converges  ($P$-a.s.)  to $z^0$.
\end{corollary}
{\bf Proof.} See Remark \ref{PoD} above.
\begin{remark}\label{FTC}{\rm
The rest of this section is concerned with the derivation of sufficient conditions to establish rate of convergence. In most applications, checking conditions of Lemma \ref{EGV} and Corollary \ref{EGV1D} below is difficult without establishing the convergence of $Z_t$ first. Therefore, although formally not required, we can assume that $Z_t\longrightarrow z^0$ convergence has already been established (using the lemmas and corollaries above or otherwise). Under this assumption, conditions for the rate of convergence below can be regarded as local in $z^0$, that is, they can be derived using certain continuity and differentiability assumptions of the corresponding functions at point $z^0$ (see examples in Section \ref{SpME}).
}\end{remark}
\begin{lemma}\label{EGV} Suppose that  $Z_t$  is a  process  defined by \eqref{TSA}. 
Let $\{C_t\}$ be a predictable positive definite $m\times m$ matrix process, and $\lambda_t^{max}$ and $\lambda_t^{min}$ be the largest and the smallest eigenvalues of $C_t$ respectively. Denote $\Delta_{t}=Z_t-z^0$. Suppose also that (V1) of Lemma \ref{RC} holds and
\begin{description}

\item[(R1)] there exists a predictable non-negative scalar process ${\cal P}_t$ such that
$$
\frac{2\Delta_{t-1}^TC_t{\gm}_t(z^0+\Delta_{t-1}) R_t(z^0+\Delta_{t-1})} {\lambda_t^{max}}+{\cal P}_t\leq -\rho_t\|\Delta_{t-1}\|^2,
$$ 
eventually, where $\rho_t$ is a predictable non-negative scalar process satisfying 
$$
\sum_{t=1}^\infty \left[ \frac {\lambda_t^{max}-\lambda_{t-1}^{min}}{\lambda_{t-1}^{min}} - \frac{\lambda_t^{max}}{\lambda_{t-1}^{min}}\rho_t\right]^+<\infty;
$$

\item[(R2)]
\end{description}
$$
\sum_{t=1}^\infty \frac{\lambda_t^{max}\left[E\left\{\Big{\|}\gm_t(z^0+\Delta_{t-1})\Big[R_t(z^0+\Delta_{t-1})+{\ve}_t(z^0+\Delta_{t-1})\Big]\Big{\|}^2 \mid{\cf}_{t-1}\right\}-{\cal P}_t\right]^+}{1+\lambda_{t-1}^{min}\|\Delta_{t-1}\|^2}<\infty.
$$
\\
\\
Then $(Z_t-z^0)^T C_t (Z_t-z^0)$ converges to a finite limit (P-a.s.).
\end{lemma}
{\bf Proof.} Let us check the conditions of Lemma \ref{RC} with $V_t(u)=u^TC_tu$. Condition (V1) is satisfied automatically.

Denote $R_t=R_t(z^0+\Delta_{t-1})$, $\gm_t=\gm_t(z^0+\Delta_{t-1})$ and $\ve_t=\ve_t(z^0+\Delta_{t-1})$. Since $V_t'(u)=2u^T C_t$ and $V_t''(u)=2C_t$, we have
$$
{\ck}_t(\Delta_{t-1})= \Delta V_t(\Delta_{t-1})+2\Delta_{t-1}^TC_t \gm_t R_t+E
\left\{[\gm_t(R_t+{\ve}_t)]^T C_t \gm_t (R_t+{\ve}_t) \mid
{\cf}_{t-1}\right\}
$$
Since $C_t$ is positive definite, $\lambda_t^{min}\|u\|^2\leq u^TC_tu\leq\lambda_t^{max}\|u\|^2$ for any $u \in {\mathbb R}^m$. Therefore
$$
 \Delta V_t(\Delta_{t-1})\leq (\lambda_t^{max}-\lambda_{t-1}^{min})\|\Delta_{t-1}\|^2.
$$
Denote
$$
\tilde{\cal P}_t=\lambda_t^{max}({\cal D}_t-{\cal P}_t)
$$
where
$$
{\cal D}_t=E\left\{\|\gm_t(R_t+{\ve}_t)\|^2 \mid{\cf}_{t-1}\right\}.
$$
Then
\begin{eqnarray}
{\ck}_t(\Delta_{t-1})
&\leq&(\lambda_t^{max}-\lambda_{t-1}^{min})\|\Delta_{t-1}\|^2+2\Delta_{t-1}^T C_t {\gm}_t R_t + \lambda_t^{max}{\cal D}_t \nonumber\\
&=&(\lambda_t^{max}-\lambda_{t-1}^{min})\|\Delta_{t-1}\|^2+2\Delta_{t-1}^T C_t {\gm}_t R_t + \lambda_t^{max}{\cal P}_t + \tilde{\cal P}_t\;. \nonumber
\end{eqnarray}
\\
By (R1), we have
$$
2\Delta_{t-1}^TC_t{\gm}_t R_t\leq - {\lambda_t^{max}}(\rho_t\|\Delta_{t-1}\|^2+{\cal P}_t).
$$
Therefore,
\begin{eqnarray}
{\cal K}_t(\Delta_{t-1})&\leq&(\lambda_t^{max}-\lambda_{t-1}^{min})\|\Delta_{t-1}\|^2 - {\lambda_t^{max}}(\rho_t\|\Delta_{t-1}\|^2+{\cal P}_t) + \lambda_t^{max}{\cal P}_t + \tilde{\cal P}_t \nonumber\\
&\leq&(\lambda_t^{max}-\lambda_{t-1}^{min}-\lambda_t^{max}\rho_t)\|\Delta_{t-1}\|^2+\tilde{\cal P}_t
=r_t\lambda_{t-1}^{min}\|\Delta_{t-1}\|^2+\tilde{\cal P}_t,\nonumber
\end{eqnarray}
where 
$$
r_t=(\lambda_t^{max}-\lambda_{t-1}^{min}-\lambda_t^{max}\rho_t)/\lambda_{t-1}^{min}.
$$
Since $\lambda_{t-1}^{min}\geq 0$, using the inequality $[a+b]^+\leq [a]^++[b]^+$, we have
$$
[{\cal K}_t(\Delta_{t-1})]^+\leq \lambda_{t-1}^{min}\|\Delta_{t-1}\|^2[r_t]^++[\tilde{\cal P}_t]^+.
$$
Also, since $V_{t-1}(\Delta_{t-1})=\Delta_{t-1}^TC_{t-1}\Delta_{t-1}\geq \lambda_{t-1}^{min}\|\Delta_{t-1}\|^2$,
\begin{eqnarray*}
\frac{[{\cal K}_t(\Delta_{t-1})]^+}{1+V_{t-1}(\Delta_{t-1})}
&\leq& \frac{[{\cal K}_t(\Delta_{t-1})]^+}{1+\lambda_{t-1}^{min}\|\Delta_{t-1}\|^2} 
\leq\frac{\lambda_{t-1}^{min}\|\Delta_{t-1}\|^2[r_t]^+}{1+\lambda_{t-1}^{min}\|\Delta_{t-1}\|^2}+\frac{[\tilde{\cal P}_t]^+}{1+\lambda_{t-1}^{min}\|\Delta_{t-1}\|^2}\\
&\leq&[r_t]^++\frac{[\tilde{\cal P}_t]^+}{1+\lambda_{t-1}^{min}\|\Delta_{t-1}\|^2}.
\end{eqnarray*}
By (R2), $\sum_{t=1}^\infty[\tilde{\cal P}_t]^+/(1+\lambda_{t-1}^{min}\|\Delta_{t-1}\|^2)<\infty$ and according to (R1)
$$
\sum_{t=1}^\infty [r_t]^+=\sum_{t=1}^\infty \left[ \frac {\lambda_t^{max}-\lambda_{t-1}^{min}}{\lambda_{t-1}^{min}} - \frac {\lambda_t^{max}}{\lambda_{t-1}^{min}}\rho_t\right]^+<\infty.
$$
Thus, 
$$
\sum_{t-1}^\infty\frac{[{\cal K}_t(\Delta_{t-1})]^+}{1+V_{t-1}(\Delta_{t-1})}<\infty,
$$ 
implying that Condition (V2)  of Lemma \ref{RC} holds. Thus, $(Z_t-z^0)^T C_t (Z_t-z^0)$ converges to a finite limit almost surely.\hfill $\blacksquare$
\begin{remark}\label{FTC2}{\rm
The choice ${\cal P}_t=0$ means that (R2) becomes more restrictive imposing stronger probabilistic restrictions on the model. Now, if $\Delta_{t-1}^TC_t{\gm}_t(z^0+\Delta_{t-1}) R_t(z^0+\Delta_{t-1})$ is eventually negative with a large absolute value, then it is possible to introduce a non-zero ${\cal P}_t$ without strengthening condition (R1). One possibility might be ${\cal P}_t=\|\gamma_tR_t\|^2$. In that case, since $\gamma_t$ and $R_t$ are predictable processes, and sequence $\ve_t$ is a martingale-difference, 
$$
E\{\|{\gm}_t(R_t+{\ve}_t)\|^2|{\cal F}_{t-1}\}=\|{\gm}_t R_t\|^2+E\{\|{\gm}_t{\ve}_t\|^2|{\cal F}_{t-1}\}.
$$
Then condition (R2) can be rewritten as
$$
\sum_{t=1}^\infty \lambda_t^{max}E\{\|{\gm}_t(z^0+\Delta_{t-1}){\ve}_t(z^0+\Delta_{t-1})\|^2|{\cal F}_{t-1}\}<\infty.
$$
}\end{remark}
\begin{remark}\label{FTC3}{\rm
The next corollary is a special case of Lemma \ref{EGV} when the step-size sequence is a sequence of scalar matrices, i.e. $\gamma_t(Z_{t-1})=a_t^{-1}\bf I$, where $a_t$ is non-decreasing and positive. 
}\end{remark}             
\begin{corollary}\label{EGV1D} Let $Z_t$ be a  process defined by \eqref{TSA}. Suppose that $a_t>0$ is a non-decreasing sequence and 
\begin{description}

\item[(W1)] 
$$
\Delta_{t-1}^T R_t(Z_{t-1}) \leq -{\frac 1 2} \Delta a_t\|\Delta_{t-1}\|^2
$$ 
eventually;

\item[(W2)] there exist $0<\delta\leq1$ such that,
$$
\sum_{t=1}^\infty a_t^{\delta-2}E\left\{\|(R_t(Z_{t-1})+{\ve}_t(Z_{t-1}))\|^2 \mid{\cf}_{t-1}\right\}<\infty.
$$
\end{description}
Then $a_t^{\delta}\|Z_t-z^0\|^2$ converges to a finite limit ($P$-a.s.).
\end{corollary}
{\bf Proof.} Consider Lemma \ref{EGV} with $\gm_t=\gm_t(z)=a_t^{-1}\bf I$, $C_t=a_t^{\delta}\bf I$, ${\cal P}_t=0$ and $\rho_t=\Delta a_t/a_{t}$. To check (R2), denote the infinite sum in (R2) by $Q$, then 
\begin{eqnarray*}
Q
&\leq& \sum_{t=1}^\infty \lambda_t^{max}\left[E\left\{\Big{\|}\gm_t\Big[R_t(z^0+\Delta_{t-1})+{\ve}_t(z^0+\Delta_{t-1})\Big]\Big{\|}^2 \mid{\cf}_{t-1}\right\}-{\cal P}_t\right]^+\\
&\leq& \sum_{t=1}^\infty \lambda_t^{max} \|\gm_t\|^2E\left\{\|(R_t(Z_{t-1})+{\ve}_t(Z_{t-1}))\|^2 \mid{\cf}_{t-1}\right\}.
\end{eqnarray*}
Now, since $\lambda_t^{min}=\lambda_t^{max}=a_t^\delta$ and $\|\gamma_t\|^2=a_t^{-2}$, condition (W2) leads to (R2).

Since $\rho_t=\Delta a_t/a_{t}<1$ and $(a_{t}/a_{t-1})^{\delta}\leq a_{t}/a_{t-1}$,
\begin{eqnarray*}
\sum_{t=1}^\infty \left[ \frac {\lambda_t^{max}-\lambda_{t-1}^{min}}{\lambda_{t-1}^{min}} - \frac{\lambda_t^{max}}{\lambda_{t-1}^{min}}\rho_t\right]^+&=&\sum_{t=1}^\infty \left[ \frac {a_t^{\delta}-a_{t-1}^{\delta}}{a_{t-1}^{\delta}} - \frac{a_t^{\delta}}{a_{t-1}^{\delta}}\rho_t\right]^+\\
&=&\sum_{t=1}^\infty \left[ (1-\rho_t) \frac{a_t^{\delta}}{a_{t-1}^{\delta}}-1\right]^+\\
&\leq&\sum_{t=1}^\infty \left[ (1-\frac{\Delta a_t}{a_{t}}) \frac{a_{t}}{a_{t-1}}-1\right]^+
=0\;.
\end{eqnarray*}
Therefore, (W1) leads to (R1). According to Remark \ref{tsphere}, condition (V1) holds since $V_t(u)=a_t^\delta\|u\|^2$. Thus, all the conditions of Lemma \ref{EGV} hold and $a_t^{\delta}\|Z_t-z^0\|^2$ converges to a finite limit ($P$-a.s.). \hfill$\blacksquare$

\begin{corollary}\label{Dev} Let $Z_t$ be a process defined by \eqref{TSA} where $z^0\in \mathbb R$, $\gamma_t(Z_{t-1})=1/t$ and the truncation sequence $U_t$ is admissible. Suppose that $Z_t\longrightarrow z^0$ and
\begin{description}
\item[(Y1)]  $R_t'(z^0)\leq-1/2$ eventually;
\item[(Y2)] $R_t(z)$ and $\sigma_t^2(z)=E(\ve_t^2(z)|{\cal F}_{t-1})$ are locally uniformly bounded at $z^0$ w.r.t. $t$; that is, there exists a constant $K$ such that $|R_t(\xi_t)|\leq K$ and $|\sigma_t^2(\xi_t)|\leq K$ eventually, for any $\xi_t\longrightarrow z^0$.
\end{description}
Then $t^{\delta}(Z_t-z^0)^2$ converges to a finite limit ($P$-a.s.), for any $\delta<1$.
\end{corollary}
{\bf Proof.} Consider Corollary \ref{EGV1D} with $a_t=t$. In the one-dimensional case, condition (W1) can be rewritten as
$$
\frac {R_t(z^0+\Delta_{t-1})}{\Delta_{t-1}} \leq -\frac{1}{2}.
$$  
Condition (W1) now follows from (Y1).

Since $E\{\ve_t(z)|{\cal F}_{t-1}\}=0$, using (Y2) we have for any $\delta<1$,
\begin{eqnarray*}
&&\sum_{t=1}^\infty t^{\delta-2}E\left\{(R_t(Z_{t-1})+{\ve}_t(Z_{t-1}))^2 \mid{\cf}_{t-1}\right\}\\
&=&\sum_{t=1}^\infty t^{\delta-2}R_t^2(Z_{t-1})+\sum_{t=1}^\infty t^{\delta-2}E\left\{{\ve}_t^2(Z_{t-1}) \mid{\cf}_{t-1}\right\}
<\infty.
\end{eqnarray*}
Thus, condition (W2) holds. Therefore, $t^{\delta}(Z_t-z^0)^2$ converges to a finite limit ($P$-a.s.), for any $\delta<1$.\hfill$\blacksquare$
\begin{remark}\label{FTC4}{\rm
Corollary \ref{Dev} gives simple but more restrictive sufficient conditions  to derive the rate of convergence in one-dimensional cases. It is easy to see that all conditions of Corollary \ref{Dev} trivially hold, if e.g., $\ve_t$ are state independent i.i.d. random variables with a finite second moment, $R_t(z)=R(z)$, and $R'(z^0)\leq-1/2$.
}\end{remark}

%
%

\section{Classical problem stochastic approximation}\label{CSA}

Consider the classical problem of stochastic approximation to find a root $z^0$ of the equation $R(z^0)=0$. Let us take a step-size sequence $\gamma_t= a_t^{-1}\bf I$, where $a_t\longrightarrow\infty$ is a predictable scalar process, and consider the procedure 
\begin{equation}\label{SN}
Z_{t}=\Phi_{U_t}\Big(Z_{t-1}+a_t^{-1} [ R(Z_{t-1})+\ve_t(Z_{t-1})]\Big).
\end{equation}
\begin{corollary} \label{ClC}
Suppose that $Z_t$  is a  process  defined by \eqref{SN}, truncation sequence $U_t$ is admissible, and
\begin{description}
\item{(H1)} 
$$
 (z-z^0)^T  R(z)\leq0
$$
for any $z\in \mathbb R^m$ with the property that $z\in U_t$ eventually;
\item{(H2)} there exists a predictable process $r_t$ such that
$$
\sup_{z\in U_{t-1}}{\| R(z)\|}\leq r_t 
\;\;\;\mbox{ where }\;\;\;
\sum_{t=1}^{\infty}a_t^{-2}r_t<\infty;
$$
\item{(H3)} there exists a predictable process $e_t$ such that
$$
\sup_{z\in U_{t-1}}\frac{E\{\|\ve_t(z)\|^2|{\cal F}_{t-1}\}}{1+\|z-z^0\|^2}\leq e_t 
$$
eventually, where
$$
\sum_{t=1}^{\infty}e_ta_t^{-2}<\infty\;\;P\mbox{-a.s.}.
$$
\end{description}
Then $\|Z_t-z^0\|$ converges to a finite limit (P-a.s.) for any initial value $Z_0$.
\\
\\
Furthermore, suppose that
\begin{description}
\item{(H4)} $R(z)$ is continuous at $z^0$ and
$
 (z-z^0)^T  R(z)<0
$
for all $z$ with the property that $z \in  U_t \backslash \{z^0\}$ eventually;
\item{(H5)} 
$$
\sum_{t=1}^\infty a_t^{-1}=\infty.
$$
\end{description}
Then $Z_t \longrightarrow z^0$ (P-a.s.).
\end{corollary}
{\bf Proof.} Consider Corollary \ref{SC} with $R_t=R$. Condition (D1) trivially holds. Since $E\left\{ \ve_t(u)\mid {\cf}_{t-1} \right\}=0$, we have 
\begin{equation*}
E\left\{ \| R(z)+\ve(z) \|^2 \mid {\cf}_{t-1} \right\} =
\|R(z)\|^2 +E\left\{  \|\ve_t(z)\|^2  \mid {\cf}_{t-1} \right\}. 
\end{equation*}
Now condition (D2) holds with $p_t=r_t+e_t$.

By (H4), there exists a constant $\nu>0$ such that for each $\epsilon \in (0,1)$
$$
\inf_{\stackrel{ \ve \le \|z-z^o\| \le 1/\ve}{z\in U_{t-1}}} -(z-z^0)^T R(u)>\nu
$$
eventually and by (H5)
$
\sum_{t=1}^\infty \nu a_t^{-1}= \nu\sum_{t=1}^\infty a_t^{-1}=\infty.
$
This implies that (D3) also holds. Therefore, by Corollary \ref{SC}, $Z_t\longrightarrow z^0$ almost surely.\hfill$\blacksquare$
\begin{remark}\label{FH3}{\rm
Suppose that   $\ve_t=\ve_t(z)$ is an  error term which does not depend on $z$ and denote 
  $$
\sigma_t^2= {E \left\{ \|\ve_t\|^2 \mid{\cf}_{t-1}\right\}}
$$
Then condition (H3) holds if 
\begin{equation}\label{statefree}
  \sum_{t=1}^{\infty} \sigma_t^2 {a_t^{-2}} <\infty, \qquad P\mbox{-a.s.}.
\end{equation}
This shows that the requirement on the error terms are quite weak. In particular, the conditional variances  do not have to be bounded w.r.t. t. 
}\end{remark}
\begin{remark} {\rm (a) If the truncation sets are uniformly bounded, then some of the conditions above can be weakened considerably. For example, condition (H2) in Corollary \ref{ClC} will automatically hold given that $\sum_{t=1}^\infty a_t^{-2}<\infty$.
\\
(b) Also if it is only required that $Z_t$ converges to any finite limit, the step-size sequence $a_t$ can go to infinity at any rate as long as $\sum_{t=1}^\infty a_t^{-2}<\infty$. However, in order to have $Z_t\longrightarrow z^0$, one must ensure that $a_t$ does not increase too fast. Also, the variances of the error terms can go to infinity as $t$ tends to infinity, as long as the sum in (H3) is bounded.
}\end{remark}
\begin{remark}\label{K-C}{\rm
To compare the above result to that of Kushner-Clark's setting, let us assume boundedness of $Z_t$. Then there exists a  compact set  $U$ such  that $Z_t \in U$. Without lost of generality, we can assume that $z^0 \in U$. Then $Z_t$  in  Corollary \ref{ClC}  can be assumed to be generated using the truncations on $U_t \cap U$. Let us assume that $\sum_{s=1}^\infty a_t^{-2} < \infty$. Then,  condition (H2) will hold if, e.g., $R(z)$ is a continuous function.  Also, in this case, given that the error terms $\varepsilon_t(z)$  are continuous in $z$ with some uniformity  w.r.t. t,  they will in fact behave in the same way as  state independent error terms. Therefore, a condition of the type \eqref{statefree} given in Remark \ref{FH3} will be sufficient for  (H3). 
}\end{remark} 
\begin{corollary} \label{ClaRC}
Suppose that $Z_t$, defined by \eqref{SN}, converges to $z^0$ (P-a.s.) and truncation sequence  $U_t$ is admissible. Suppose also that
\begin{description}
\item{(B1)} 
$$
u^T R(z^0+u)\leq- \frac 1 2 \|u\|^2\;\;\;\mbox{ for small $u$'s;}
$$
\item{(B2)} $a_t>0$ is non-decreasing with
$$
\sum_{t=1}^{\infty}\left[\frac{\Delta a_{t}-1}{a_{t-1}}\right]^+<\infty ;
$$
\item{(B3)} there exist $\delta\in(0,1)$ such that 
$$
\sum_{t=1}^\infty a_t^{\delta-2}\|R(z^0+v_t)\|^2<\infty ~~~~~
\mbox{and} ~~~~~
\sum_{t=1}^\infty a_t^{\delta-2}E\{\|\ve_t(z^0+v_t)\|^2|{\cal F}_{t-1}\}<\infty,
$$
where $v_t\in U_t$ is any predictable process with the property $v_t\longrightarrow 0$.
\end{description}
Then $a_t^{\delta} \|Z_t-z^0\|^2$ converges (P-a.s.) to a finite limit.
\end{corollary} 
{\bf Proof.} Let us check that conditions of Lemma \ref{EGV} hold with $R_t=R$, $\rho_t=a_t^{-1}$, ${\cal P}_t=0$ and $C_t=a_t^{\delta} {\bf I}$. We have $\lambda_t^{max}=\lambda_t^{min}=a_t^{\delta}$ by (B2), and
\begin{eqnarray*}
&&\sum_{t=1}^\infty\left[\frac{\lambda_t^{max}-\lambda_{t-1}^{min}}{\lambda_{t-1}^{min}}-\frac{\lambda_t^{max}}{\lambda_{t-1}^{min}}\rho_t\right]^+ = \sum_{t=1}^\infty\left[\frac{a_t^{\delta}-a_{t-1}^{\delta}}{a_{t-1}^{\delta}}-\frac{a_t^{\delta}}{a_{t-1}^{\delta}a_t}\right]^+\\
&=&\sum_{t=1}^\infty\left[\left(\frac{a_t}{a_{t-1}}\right)^{\delta}(1-a_t^{-1})-1\right]^+
\leq\sum_{t=1}^\infty\left[\frac{a_{t}}{a_{t-1}}(1-a_{t}^{-1})-1\right]^++C\\
&=&\sum_{t=1}^\infty\left[\frac{\Delta a_{t}-1}{a_{t-1}}\right]^++C
<\infty
\end{eqnarray*}
for some constant C. So (B1) leads to (R1). Also since $Z_t \longrightarrow z^0$,
\begin{eqnarray*}
&&\sum_{t=1}^\infty \frac{\lambda_t^{max}\left[E\left\{\|\gm_t(R_t+{\ve}_t)\|^2 \mid{\cf}_{t-1}\right\}-{\cal P}_t\right]^+}{1+\lambda_{t-1}^{min}\|\Delta_{t-1}\|^2}\\
&\leq&\sum_{t=1}^\infty \lambda_t^{max}\left[E\left\{\|\gm_t(R_t+{\ve}_t)\|^2 \mid{\cf}_{t-1}\right\}-{\cal P}_t\right]^+\\
&=&\sum_{t=1}^\infty a_t^{\delta}E\left\{\|a_t^{-1}(R_t+{\ve}_t)\|^2 \mid{\cf}_{t-1}\right\}\\
&\leq&\sum_{t=1}^\infty a_t^{\delta-2}\| R(Z_{t-1})\|^2+\sum_{t=1}^\infty a_t^{\delta-2}E\{\|{\ve_t(Z_{t-1})}\|^2|{\cal F}_{t-1}\},
\end{eqnarray*}
condition (R2) follows from (B3). Therefore by Lemma \ref{EGV}, $(Z_t-z^0)^TC_t(Z_t-z^0)=a_t^\delta\|Z_t-z^0\|\longrightarrow0$ ($P$-a.s.).\hfill $\blacksquare$
\begin{remark}\label{pow}{\rm
It follows from Proposition \ref{ADT} in Appendix  that if $a_t=t^\epsilon$ with $\epsilon>1$, then (B2) doesn't hold. However, condition (B2) holds if $a_t=t^\epsilon$ for all $\epsilon\leq 1$. Indeed, 
\begin{eqnarray*}
\sum_{t=1}^{\infty}\left[\frac{\Delta a_{t}-1}{a_{t-1}}\right]^+&=&\sum_{t=1}^{\infty}\left[\left(\frac{t}{t-1}\right)^\epsilon-1-\frac1{(t-1)^\epsilon}\right]^+\\
&\leq&\sum_{t=1}^{\infty}\left[\frac{t}{t-1}-1-\frac1{t-1}\right]^+
=0.
\end{eqnarray*}
}\end{remark}


\begin{corollary}\label{Po}
Suppose that $Z_t\longrightarrow z^0$, where $Z_t$ is defined by \eqref{SN} with $a_t=t^\epsilon$ where $\epsilon\in(1/2,1]$, and (B1) in Corollary \ref{ClaRC} holds. Suppose also that $R$ is continuous at $z^0$ and there exists $0<\delta <2-1/\epsilon$ such that
\begin{description}
\item[(BB)]
$$
\sum_{t=1}^\infty \frac{1}{t^{(2-\delta)\epsilon}}E\{\|\ve_t(z^0+v_t)\|^2|{\cal F}_{t-1}\}<\infty.
$$
where $v_t\in U_t$ is any predictable process with the property $v_t\longrightarrow 0$.
\end{description}
Then $t^{\delta}\|Z_t-z^0\|^2$ converges to a finite limit (P-a.s.).
\end{corollary}
{\bf Proof.} Let us check conditions of Corollary \ref{ClaRC} with $a_t=t^\epsilon$ where $\epsilon\in(1/2,1]$. Condition (B2) is satisfied (See Remark \ref{pow}). Since $R$ is continuous at $z^0$ and $Z_t \longrightarrow z^0$, it follows that $R(z^0+v_t)$ in (B3) is bounded. Also, $a_t^{\delta-2}=t^{(\delta-2)\ve}$ and since $(\delta-2)\epsilon<-1$, it follows that the first part of (B3) holds. The second part is a consequence of (BB). The result is now immediate from Corollary \ref{ClaRC}.\hfill$\blacksquare$
\begin{remark}\label{Po1}{\rm
Suppose that $a_t=t^{\ve}$ with $\ve\in(1/2,1)$ and $\sup_t E\{\|\ve_t(z)\|^2|{\cal F}_{t-1}\}<\infty$ (e.g., assume that $\ve_t=\ve_t(z)$ are state independent and i.i.d.). Then, since $(\delta-2)\epsilon<-1$, condition (BB) in Corollary \ref{Po} automatically holds for any $\delta<2-1/\epsilon$. It therefore follows that the step-size sequence $a_t=t^\epsilon$, $\epsilon\in(1/2,1)$ produces SA procedures which converge with the rate $t^{-\alpha}$ where $\alpha<1-\frac1 {2\epsilon}$. For example, the step-size $a_t=t^{3/4}$ would produce the SA procedures, which converge with the rate $t^{-1/3}$.
}\end{remark}

%
%

\section{Special models and examples }\label{SpME}

\subsection{Finding a root of a polynomial}\label{Poly}
Let  $l$  be a positive integer and 
 $$
 R(z)=-\sum_{i=1}^{l}C_i(z-z^0)^i,
 $$                
where $z , z^0 \in \mathbb{R}$ and $C_i$ are real constants. Suppose that
$$
(z-z^0)R(z)\leq0 \;\;\;\mbox{ for all } \;\;\;z\in \mathbb R.
$$
Note that if $l >1$, the SA without truncations fails to satisfy the standard condition on the rate of growth at infinity. Therefore, one needs to use  slowly expanding truncations to  slow down the growth of $R$ at infinity.
Consider $Z_t$ defined by \eqref{SN} with a truncation sequence $U_t=[-u_t, u_t]$, where $u_t \longrightarrow \infty$ is a sequence of non-decreasing positive numbers.   Suppose that  
\begin{equation}\label{ut}
    \sum_{t=1}^{\infty}  u_t^{2l}~ a_t^{-2} <\infty.
\end{equation}
 Then, provided that the measurement errors satisfy condition (H3)  of Corollary \ref{ClC},  $|Z_t-z^0|$ converges ($P$-a.s.) to a finite limit.
 
 Indeed,  condition (H1) of Corollary \ref{ClC} trivially holds. For large $t$'s,
\begin{eqnarray*} 
&&\sup_{z \in [-u_{t-1}, u_{t-1}]} ~ \|R(z)  \|^2\;\;\leq  
\sup_{z \in [-u_{t-1}, u_{t-1}]} \left[\sum_{i=1}^l C_i(z-z^0)^{i}\right]^2   \\
&\leq& \sup_{z \in [-u_{t-1}, u_{t-1}]} \sum_{i=1}^l C_i^2(z-z^0)^{2i}  
\leq  \sum_{i=1}^l C_i^2(2u_{t-1})^{2i}   
\leq l4^l C_l^2 u_{t-1}^{2l},
\end{eqnarray*}
which, by \eqref{ut}, implies  condition (H2)  of Corollary \ref{ClC}.

Furthermore, if $z^0$ is a unique root, then provided that 
\begin{equation}\label{SoA}
\sum_{t=1}^{\infty}   a_{t}^{-1}=\infty,
\end{equation}
it follows from Corollary \ref{ClC} that $Z_t\longrightarrow z^0$ ($P$-a.s.).
One can always choose a suitable truncation sequence which satisfies \eqref{ut} and \eqref{SoA}. For example, if the degree of the polynomial is known to be $l$ (or at most $l$), and $a_t=t$,  then one can take $u_t=Ct^{r/{2l}}$, where $C$ and $r$ are some positive constants and $r < 1$.  One can also take a truncation sequence which is independent of $l$, e.g.,  $u_t=C \log t$, where $C$ is a positive constant.

Suppose also that 
$$
C_1\geq \frac1 2,\;\;\;a_t=t^{\epsilon}\;\;\;\mbox{where} \;\;\;\epsilon\in(0,1]
$$
and condition (BB) in Corollary \ref{Po} holds (e.g., one can assume for simplicity that $\ve_t$'s are state independent and i.i.d.). Then $t^\alpha (Z_t-z^0)\xrightarrow{a.s.} 0$ for any $\alpha<1-1/2\epsilon$.

Indeed, since $R'(z^0)=-C_1\leq  -1/2$, condition (B1) of Corollary \ref{ClaRC} holds. Now, the above convergence is a consequence of Corollary \ref{Po} and Remark \ref{Po1}.


\subsection{Linear procedures}
Consider the recursive procedure
\begin{equation}\label{LP}
Z_t=Z_{t-1}+\gamma_t(h_t-\beta_t Z_{t-1})
\end{equation}
where $\gamma_t$ is a predictable positive definite matrix process, $\beta_t$ is a predictable positive semi-definite matrix process and $h_t$ is an adapted vector process (i.e., $h_t$ is ${\cal F}_t$-measurable for $t\geq 1$). If we assume that $E\{h_t|{\cal F}_{t-1}\}=\beta_t z^0$, we can view \eqref{LP} as a SA procedure designed to find the common root $z^0$ of the linear functions
$$
R_t(u)=E\{h_t-\beta_t u|{\cal F}_{t-1}\}=E\{h_t|{\cal F}_{t-1}\}-\beta_t u=\beta_t(z^0-u)
$$ 
which is observed with the random noise 
$$\ve_t(u)=h_t-\beta_t u-R_t(u)=h_t-E\{h_t|{\cal F}_{t-1}\}=h_t-\beta_t z^0.$$

\begin{corollary}\label{LRA}
Suppose that $Z_t$ is defined by \eqref{LP} with $E(h_t|{\cal F}_{t-1})=\beta_tz^0$. Suppose also that  $a_t$ is a non-decreasing positive predictable process and

\begin{description}
\item[(G1)]
$\Delta\gamma_t^{-1}-2\beta_t+\beta_t\gamma_t\beta_t$ is negative semi-definite eventually;
\item[(G2)] 
$$
\sum_{t=1}^{\infty} a_t^{-1}E\{(h_t-\beta_t z^0)^T\gamma_t(h_t-\beta_t z^0)|{\cal F}_{t-1}\}<\infty.
$$

\end{description}
Then $a_t^{-1}(Z_t-z^0)^T\gamma_t^{-1}(Z_t-z^0)$ converges to a finite limit (P-a.s.). 
\end{corollary}
{\bf Proof.}
Let us show that conditions of Lemma \ref{RC} hold with $V_t(u)=a_t^{-1}u^T\gamma_t^{-1}u$. Condition (V1) trivially holds. We have $V'_t(u)=2a_t^{-1}u^T\gamma_t^{-1}$, $V''_t(u)=2a_t^{-1}\gamma_t^{-1}$, $R_t(z^0+u)=-\beta_tu$ and $R_t(u)+\ve_t(u)=h_t-\beta_tu$. Since $E(h_t-\beta_tz^0|{\cal F}_{t-1})=0$, for $\eta_t$ defined in (V2) we have
\begin{eqnarray*}
\eta_t(z^0+u)&=& a_t^{-1}E\Big\{(h_t-\beta_tz^0-\beta_tu)^T\gamma_t(h_t-\beta_tz^0-\beta_tu)\Big|{\cal F}_{t-1}\Big\}\\
&=& a_t^{-1}E\Big\{(h_t-\beta_tz^0)^T\gamma_t(h_t-\beta_tz^0)\Big|{\cal F}_{t-1}\Big\}+a_t^{-1}(\beta_tu)^T\gamma_t(\beta_tu)\;.
\end{eqnarray*}
Also, 
$$
\Delta V_t(u)=u^T[(a_t\gamma_t)^{-1}-(a_{t-1}\gamma_{t-1})^{-1}]u\leq u^T(a_t\gamma_t)^{-1}u-u^T(a_t\gamma_{t-1})^{-1}u=u^Ta_t^{-1}\Delta\gamma_t^{-1}u.
$$
Denoting
$$
{\cal J}_t=a_t^{-1}E\Big\{(h_t-\beta_tz^0)^T\gamma_t(h_t-\beta_tz^0)\Big|{\cal F}_{t-1}\Big\},
$$
for ${\cal K}_t$ from (V2), we have
\begin{eqnarray*}
{\cal K}_t(u)&\leq& a_t^{-1}u^T\Delta\gamma_t^{-1}u-2 a_t^{-1}u^T\beta_tu+a_t^{-1}u^T\beta_t^T\gamma_t\beta_tu+{\cal J}_t\\
&=& a_t^{-1}u^T(\Delta\gamma_t^{-1}-2\beta_t+\beta_t^T\gamma_t\beta_t)u+{\cal J}_t\;.
\end{eqnarray*}
Condition (V2) is now immediate from (G1) and (G2) since 
$$
[1+V_{t-1}(\Delta_{t-1})]^{-1}[{\cal K}_t(\Delta_{t-1})]^+\leq[{\cal K}_t(\Delta_{t-1})]^+\leq{\cal J}_t\;.
$$
Thus, all the conditions of Lemma \ref{RC} hold which implies the required result.
\hfill$\blacksquare$

\begin{corollary}\label{LR2}
Suppose that $\Delta\gm_t^{-1}=\beta_t$, then (G1) in Corollary \ref{LRA} holds.
\end{corollary}
{\bf Proof.} Since $\Delta\gm_t^{-1}$ is positive semi-definite, it follows that $\Delta\gm_t$ is negative semi-definite \Big(see Horn and Johnson (1985)\nocite{HaJ} Corollary 7.7.4(a)\Big). Also since
\begin{eqnarray*}
\Delta\gamma_t^{-1}-2\beta_t+\beta_t\gamma_t\beta_t&=&-\Delta\gamma_t^{-1}+\Delta\gamma_t^{-1}\gamma_t\Delta\gamma_t^{-1}=-\Delta\gamma_t^{-1}+\gamma_t^{-1}-2\gamma_{t-1}^{-1}+\gamma_{t-1}^{-1}\gm_t\gm_{t-1}^{-1}\\
&=&-\gamma_{t-1}^{-1}+\gamma_{t-1}^{-1}(\gm_{t-1}+\Delta\gm_t)\gm_{t-1}^{-1}
=\gamma_{t-1}^{-1}\Delta\gm_t\gm_{t-1}^{-1}, 
\end{eqnarray*}
it follows that (G1) holds. \hfill  $ \blacksquare $


\subsection{Parameter estimation in Autoregressive models}\label{REAR}

Consider an AR(m) process
\begin{equation}\label{ARMo}
X_t=\theta^{(1)}X_{t-1}+\theta^{(2)}X_{t-2}+\dots+\theta^{(m)}X_{t-m}+\xi_t=\theta^TX_{t-m}^{t-1}+\xi_t \nonumber
\end{equation}
where $\theta=(\theta^{(1)},...,\theta^{(m)})^T$, $X_{t-m}^{t-1}=(X_{t-1},..., X_{t-m})^T$ and $\xi_t$ is a 
martingale-difference (i.e., $E\{ \xi_t |{\cal F}_{t-1} \}=0$). If the pdf of    ${\xi}_t$ w.r.t. Lebesgue's measure is $g_t(x)$, then the conditional probability density function of $X_t$ given the past observations is
$$
f_t(x,\theta|X_{1}^{t-1})=f_t(x,\theta|X_{t-m}^{t-1})=g_t(x-\theta^TX_{t-m}^{t-1})
$$
and
\begin{equation}\label{Argmle}
\frac{ {f'}_t^T(\theta, x | X_1^{t-1})}
 {f_ t(\theta, x |X_1^{t-1})}
=-\frac{g_t'(x-\theta^TX_{t-m}^{t-1})}{g_t(x-\theta^TX_{t-m}^{t-1})}X_{t-m}^{t-1}. \nonumber
\end{equation}
It is easy to see that the conditional Fisher information  \eqref{fish} is
$$
I_t=\sum_{s=1}^{t}l_{gs} X_{s-m}^{s-1}(X_{s-m}^{s-1})^T ~~~~\mbox{where} ~~~~
l_{gt}=\int_{-\infty}^{\infty} \left(\frac{g_t'(x)}{g_t(x)}\right)^2g_t(x)dx.
$$
 The inverse $ I_t^{-1}$   can also be generated recursively  by
\begin{equation}\label{RIo}
 I_t^{-1}= I_{t-1}^{-1}- l_{gt} I_{t-1}^{-1}X_{t-m}^{t-1}(1+ l_{gt}(X_{t-m}^{t-1})^T I_{t-1}^{-1}X_{t-m}^{t-1})^{-1}(X_{t-m}^{t-1})^T I_{t-1}^{-1}.
\end{equation}
(Note that this can be derived either directly,  or using the matrix inversion formula, sometimes referred to as the Sherman-Morrison formula.)

Thus, the on-line likelihood procedure in this case can be derived by the following recursion
\begin{equation}\label{ARMLo}
\hat \theta_t=\hat \theta_{t-1}- I_t^{-1}X_{t-m}^{t-1}\frac{g_t'}{g_t}(X_t- \hat \theta_{t-1}^T X_{t-m}^{t-1})
\end{equation}                       
where $ I_t^{-1}$ is also derived on-line using formula \eqref{RIo}.
In general, to include robust estimation procedures, and also to use any available auxiliary information, one can use the following class of procedures                           
\begin{equation}\label{ArgM}
\hat \theta_t=\Phi_{U_t}\left(\hat \theta_{t-1}+\gamma_t H(X_{t-m}^{t-1})
 \varphi_t(X_t-\hat \theta_{t-1}^TX_{t-m}^{t-1})\right),
\end{equation}
where
$\varphi_t:\mathbb R\mapsto\mathbb R$ and $H:\mathbb R^m\mapsto\mathbb R^m$  are suitably chosen functions and $\gm_t$ is an $m\times m$ matrix valued step-size sequence.


\begin{example}{\it (Recursive least squares procedures)}\label{RLS} 
{\rm Recursive least squares (RLS) estimator of $\theta=(\theta^{(1)},\dots,\theta^{(m)})^T$ is generated by the following procedure
\begin{equation}\label{LS}
\hat\theta_t=\hat\theta_{t-1}+\hat I_t^{-1}X_{t-m}^{t-1}[X_t-( X_{t-m}^{t-1})^T\hat\theta_{t-1}],
\end{equation}
\begin{equation}\label{RI}
\hat I_t^{-1}=\hat I_{t-1}^{-1}-\hat I_{t-1}^{-1}X_{t-m}^{t-1}[1+(X_{t-m}^{t-1})^T\hat I_{t-1}^{-1}X_{t-m}^{t-1}]^{-1}(X_{t-m}^{t-1})^T \hat I_{t-1}^{-1},
\end{equation}
where $\hat\theta_0$ and  a positive definite $\hat I_0^{-1}$ are some starting values. 
Note that \eqref{LS} is a particular case of \eqref{ArgM}, and it also coincides with the maximum likelihood procedure \eqref{ARMLo} in the case when $\xi_t$ are i.i.d. Gaussian r.v.'s.


\begin{corollary}\label{RCAR}

Consider $\hat\theta_t$ defined by \eqref{LS} and \eqref{RI}. Suppose that there exists a non-decreasing sequence $a_t>0$ such that
\\
$$
\sum_{t=1}^{\infty} a_t^{-1} (X_{t-m}^{t-1})^T \hat I_t^{-1} X_{t-m}^{t-1}E\{\xi_t^2|{\cal F}_{t-1}\}<\infty.
$$
Then $a_t^{-1} (\hat\theta_t-\theta)^T\hat I_t (\hat\theta_t-\theta)$ converges to a finite limit ($P^\theta$-a.s.).
\end{corollary}
{\bf Proof.}
Let us check the condition of Corollary \ref{LRA}. Obviously, the matrix $\gamma_t=\hat I_t^{-1}=\hat I_0^{-1}+\sum_{s=1}^t X_{s-m}^{s-1}(X_{s-m}^{s-1})^T$ is positive definite and $\Delta \hat{I}_t^{-1}=\beta_t=X_{t-m}^{t-1}(X_{t-m}^{t-1})^T$ is positive semi-definite. By Corollary \ref{LR2}, condition (G1) holds.
We also have
\begin{eqnarray*}
\sum_{t=1}^{\infty} a_t^{-1}E\{(h_t-\beta_t z^0)^T\gamma_t(h_t-\beta_t z^0)|{\cal F}_{t-1}\}
=\sum_{t=1}^{\infty} a_t^{-1}E\{\xi_t(X_{t-m}^{t-1})^T\hat I_t^{-1}X_{t-m}^{t-1}\xi_t|{\cal F}_{t-1}\}\\
=\sum_{t=1}^{\infty} a_t^{-1} (X_{t-m}^{t-1})^T\hat I_t^{-1}X_{t-m}^{t-1}E\{\xi_t^2|{\cal F}_{t-1}\}<\infty.
\end{eqnarray*}
So condition (G2) holds. Hence all conditions of Corollary \ref{LRA} hold which completes the proof.
\hfill$\blacksquare$

\begin{corollary}\label{RCAR2}
Consider $\hat\theta_t$ defined by \eqref{LS} and \eqref{RI}. Suppose that 
\begin{description}
\item[(P1)] there exists a non-decreasing sequence $\kappa_t\longrightarrow \infty$ such that 
$$
{\hat I_t}/{\kappa_t} \longrightarrow G
$$
where $G<\infty$ is a positive definite $m\times m$ matrix;
\item[(P2)] there exists $\epsilon^0\in [0,1)$ such that
$$
E\left\{ \xi_t^2|{\cal F}_{t-1}  \right\}\leq \kappa_t^{\epsilon^0}\;\;\;\mbox{ eventually.}
$$
\end{description}
Then $\kappa_t^{1-\delta}\|\hat\theta_t-\theta\|^2\longrightarrow0$ ($P^\theta$-a.s.) for all $\delta\in(\epsilon^0,1]$.
\end{corollary}
{\bf Proof.} Consider Corollary \ref{RCAR} with $a_t=\kappa_t^{\delta}$ for a certain $\delta\in(\epsilon^0,1]$. By (P2), there exists $t^0$ such that
$$
\sum_{t=t^0}^{\infty} a_t^{-1} (X_{t-m}^{t-1})^T \hat I_t^{-1} X_{t-m}^{t-1}E\{\xi_t^2|{\cal F}_{t-1}\}\leq\sum_{t=t^0}^{\infty} \kappa_t^{\epsilon^0-\delta} (X_{t-m}^{t-1})^T \hat I_t^{-1} X_{t-m}^{t-1}
$$
eventually. Now, using (P1) and Lemma \ref{KapPr} in Appendix , the above sum converges to a finite limit implying conditions of Corollary \ref{RCAR} hold. 
Therefore, $(\hat{\theta}_t-\theta)^T\hat I_t (\hat{\theta}_t-\theta)/{\kappa_t^{\delta}}$ tends to a finite limit. Now, the assertion of the corollary follows since $\hat I_t/\kappa_t$ converges to a finite matrix. 
\hfill$\blacksquare$
\begin{remark}\label{KaR}{\rm
{\bf (a)}
If $X_t$ is a strongly stationary process, condition (P1) will  trivially hold with  $\kappa_t=t$.
However,  using the results given above, convergence can be derived without the stationarity requirement as long as $\kappa_t^{-1} \sum_{t=1}^\infty X_{t-m}^{t-1} (X_{t-m}^{t-1})^T$ tends to a positive define matrix.
\\
{\bf (b)}  Condition (P2) demonstrates that the requirements on the innovations $\xi_t$ are quite week.  In particular, the conditional variances of the innovations do not have to be bounded w.r.t. $t$. For example, if   $\kappa_t=t$  and the variances go to infinity not faster than $t^{\varepsilon_0}$  (for some $0\le{\varepsilon_0}<1$), then  it follows that    $t^{1-\delta}\|\hat\theta_t-\theta\|^2 \to 0$  for any  $\delta \in ({\varepsilon_0}, 1)$.
\\
{\bf (c)} It follows from  (a) and (b) above that in the case of   a strongly stationary $X_t$ with 
iid innovations, $t^{1-\delta} \|\hat\theta_t-\theta\|^2 \to 0$  for any  $\delta>0$ without any additional assumptions.
}\end{remark}
}\end{example}



\numberwithin{equation}{section}
\numberwithin{lemma}{section}
\addcontentsline{toc}{chapter}{Appendix}

\section{Appendix}\label{App}

\begin{lemma}\label{RMTh} Let ${\cal F}_0$, ${\cal F}_1$, \dots be an non-decreasing sequence of $\sigma$-algebras and $X_n$, $\beta_n$, $\xi_n$, $\zeta_n \in {\cal F}_n$, $n\geq0$, be non-negative random valuables such that
$$
E(X_n|{\cal F}_{n-1}) \leq X_{n-1}(1+\beta_{n-1})+\xi_{n-1}-\zeta_{n-1},~~~n\geq 1
$$
eventually. Then
$$
\left\{ \sum_{i=1}^\infty \xi_{i-1} <\infty \right\}\cap \left\{\sum_{i=1}^\infty \beta_{i-1} <\infty\right\}\subseteq \left\{X \rightarrow \right\} \cap\left\{ \sum_{i=1}^\infty \zeta_{i-1} < \infty \right\} \qquad
 P\mbox{-a.s.},
$$
where $\{X \rightarrow \}$ denotes the set where $\lim_{n\rightarrow \infty} X_n$ exists and is finite.
\end{lemma}
{\bf Proof.} The proof can be found in Robbins and Siegmund (1985\nocite{Rob2}). Note also that this lemma is a special case of the theorem on the convergence sets of non-negative semi-martingales (see, e.g., Lazrieva et al (1997\nocite{laz2})).\hfill$\blacksquare$

\begin{proposition}\label{TsphPr}
Consider a closed sphere $U=S(\alpha, r)$ in $\mathbb R^m$ with the center at $\alpha\in\mathbb R^m$ and the radius $r$. Let $z^0\in U$ and $z\notin U$. Denote by $z'$ the closest point form $z$ to $U$, that is,
$$
z'=\alpha+\frac{r}{\|z-\alpha\|}(z-\alpha).
$$
Suppose also that $C$ is a positive definite matrix such that
\begin{equation}\label{EigA}
\lambda_C^{max}v^2\leq\lambda_C^{min}r^2, \nonumber
\end{equation}
where $\lambda_C^{max}$ and $\lambda_C^{min}$ are the largest and smallest eigenvalues of $C$ respectively and $v=\|\alpha-z^0\|$. Then  $$(z'-z^0)^TC(z'-z^0)\leq(z-z^0)^TC(z-z^0).$$
\end{proposition}
{\bf Proof.}
For $u,v\in \mathbb R^m$, define
$$
\|u\|_C=(u^TCu)^{1/2} \;\;\;\mbox{ and }\;\;\;(u,v)_C=(u^TCv)^{1/2}.
$$ 
We have
\begin{eqnarray}\label{mod}
\left|(z_0 - \alpha, z' - \alpha)_C\right| \le \|z_0 - \alpha\|_C \|z' - \alpha\|_C \le
 \sqrt{\lambda^{\mathrm max}_C}\, v \|z' - \alpha\|_C  \nonumber \\
 \le 
\sqrt{\lambda^{\mathrm min}_C}\, r \|z' - \alpha\|_C =
\sqrt{\lambda^{\mathrm min}_C}\, \|z' - \alpha\| \|z' - \alpha\|_C 
\le \|z' - \alpha\|_C^2\, . 
\end{eqnarray}
Since $z\notin U$, we have 
$$
z'=\alpha+\frac{r}{\|z-\alpha\|}(z-\alpha)=(1-\delta)\alpha+\delta z,
$$
where $\delta={r}/{\|z-\alpha\|}<1$. Then, since 
$$
z - z' = (1 - \delta) (z - \alpha) , \ \ \  z' - \alpha = \delta (z - \alpha),  \ \ \
z - z' = \frac{1 - \delta}{\delta}\, (z' - \alpha), 
$$
by \eqref{mod},
\begin{eqnarray}
(z' - z_0, z - z')_C = (z' - \alpha, z - z')_C + (\alpha - z_0, z - z')_C \nonumber\\
=\frac{1 - \delta}{\delta}\, \|z' - \alpha\|_C^2 - \frac{1 - \delta}{\delta}\, (z_0 - \alpha, z' - \alpha)_C \ge 0. \nonumber
\end{eqnarray}
Therefore, since
$
z' - z_0 = (z - z_0) - (z - z') ,
$
we get
\begin{eqnarray*}
\|z' - z_0\|_C^2 = \|z - z_0\|_C^2 + \|z - z'\|_C^2 - 2 (z - z_0, z - z')_C \\
= \|z - z_0\|_C^2 + \|z - z'\|_C^2 - 2\|z - z'\|_C^2 - 2 (z' - z_0, z - z')_C \\
= \|z - z_0\|_C^2  - \|z - z'\|_C^2 - 2 (z' - z_0, z - z')_C \le \|z - z_0\|_C^2. \\ \quad \blacksquare
\end{eqnarray*}

\begin{proposition}\label{ADT}
Suppose $a_t$, $t \in \mathbb{N}$ is a non-decreasing sequence of positive numbers such that
$$
\sum_{t = 1}^\infty \frac{1}{a_t} < \infty .
$$
Then
$$
\sum_{t = 1}^\infty \left[\frac{a_{t + 1} - a_t - 1}{a_t}\right]^+ = +\infty .
$$ 
\end{proposition}
{\bf Proof.}
Since
$$
\sum_{t = 1}^\infty \left[\frac{a_{t + 1} - a_t - 1}{a_t}\right]^+ \ge 
\sum_{t = 1}^\infty \frac{a_{t + 1} - a_t}{a_t} - \sum_{t = 1}^\infty \frac{1}{a_t}
$$
and the last series converges, it is sufficient to show that
$$
\sum_{t = 1}^\infty \frac{a_{t + 1} - a_t}{a_t} = +\infty .
$$
Note that for $b\geq a>0$, we have
$$
\frac{b - a}{a} = \int_a^b \frac{1}{a}\, d\tau \ge \int_a^b \frac{1}{\tau}\, d\tau = \ln b - \ln a .
$$
So,
\begin{eqnarray*}
\sum_{t = 1}^N \frac{a_{t + 1} - a_t}{a_t} \ge \sum_{t = 1}^N \left(\ln a_{t + 1} - \ln a_t\right) = \ln a_{N + 1} - \ln a_1 \to +\infty \
\mbox{ as } N \to \infty .  \quad  \blacksquare 
\end{eqnarray*}


\begin{lemma} \label{KapPr}
Suppose $\{\alpha_t\}$ is a sequence of real $m \times 1$ column vector, $I_t={\bf I}+\sum_{s=1}^t \alpha_s\alpha_s^T$ diverges and  $\kappa_t$ is a sequence of positive numbers satisfying:
$$
{I_t}/ {\kappa_t}\rightarrow G,
$$
where $G$ is a finite positive definite $m \times m$ matrix. Then  
$$
\sum_{t=N}^\infty \frac 1 {\kappa_t^{\delta}} \alpha_t^T I_t^{-1} \alpha_t <\infty
$$
for any $\delta>0$.
\end{lemma}
{\bf Proof.}  Since
$tr(I_t)=m+\sum_{s=1}^t \alpha_s^T\alpha_s$ is a non-decreasing sequence of positive numbers, we have (see Proposition A2 in Sharia (2007)\nocite{Shar2})
$$
\sum_{t=1}^\infty \frac {\alpha_t^T \alpha_t} {[tr(I_t)]^{1+\delta}}<\sum_{t=1}^\infty \frac {\alpha_t^T \alpha_t} {(\sum_{s=1}^t \alpha_s^T\alpha_s)^{1+\delta}}<\infty.
$$
Since {$  {I_t}/ {\kappa_t}$} converges, we have that {$ {tr(I_t)}/ {\kappa_t}$} tends to a finite limit, and
$$
\sum_{t=1}^\infty \frac {\alpha_t^T \alpha_t} {{\kappa_t}^{1+\delta}}=\sum_{t=1}^\infty \frac {\alpha_t^T \alpha_t} {tr(I_t)^{1+\delta}}\left[\frac {tr(I_t)} {\kappa_t}\right]^{1+\delta}<\infty
$$
Finally, since $G_t$ is positive definite and we have ${\kappa_t}I_t^{-1} \rightarrow G^{-1}$, and it follows that
$
{\kappa_t}\lambda_t^{max} 
$
converges to a finite limit,
where $\lambda_t^{max}$ is the largest eigenvalue of $I_t^{-1}$. Thus,
$$
\sum_{t=1}^\infty \frac 1 {\kappa_t^\delta} \alpha_t^T I_t^{-1} \alpha_t \leq\sum_{t=1}^\infty\frac {\alpha_t^T \alpha_t} {\kappa_t^{1+\delta}}\cdot \kappa_t\lambda_t^{max}< \infty. \quad \blacksquare 
$$



 

\bibliographystyle{acm}
\bibliography{bibLT1}

\end{document}